\documentclass[a4paper,twoside,11pt]{article}

\usepackage[francais,english]{babel}
\usepackage{amssymb}

\def \elem(#1,#2){  {{#1}\over \overline {\ #2\ }} }
\def \a{\kern+.6ex\lower.42ex\hbox{$\scriptstyle \iota$}\kern-1.20ex a}

\newcommand{\Teis}[2]{
   \setlength{\unitlength}{1ex}
   \begin{picture}(2,0)(0,0.4)
      \put(0,1.1){\line(1,0){2}}
      \put(0,0.9){\line(1,0){2}}
      \put(1,1.2){\makebox(0,0)[b]{$\scriptstyle #1$}}
      \put(1,0.8){\makebox(0,0)[t]{$\scriptstyle #2$}}
   \end{picture}}

\newcommand{\Teisssr}[4]{
   \setlength{\unitlength}{1ex}
   \begin{picture}(#3,3)(0,0.4)
      \put(0,1.15){\line(1,0){#3}}
      \put(0,0.85){\line(1,0){#3}}
      \put(#4,1.3){\makebox(0,0)[b]{$#1$}}
      \put(#4,0.7){\makebox(0,0)[t]{$#2$}}
   \end{picture}}

\newcommand{\bDelta}{{\Delta}}

\newcommand{\bN}{{\mathbb N}}

\newcommand{\bR}{{\mathbb R}}
\newcommand{\bC}{{\mathbb C}}
\newcommand{\cN}{{\cal N}}

\newtheorem{teorema}{Theorem}[section]
\newtheorem{lema}[teorema]{Lemma}
\newtheorem{coro}[teorema]{Corollary}
\newtheorem{prop}[teorema]{Proposition}
\newtheorem{nota}[teorema]{Remark}
\newtheorem{defi}[teorema]{Definition}

\title{Characterization of non-degenerate plane curve singularities
\footnotetext{
\begin{minipage}[t]{4.5in}
2000 MS Classification: 32S55, 14H20.\\
Key words: {\em non-degenerate plane curve singularities, Milnor number,
Newton number}\\
This research was partially supported by Spanish Projet MEC
PNMTM2004-00958
\end{minipage}}}

\author{Evelia R. Garc\'{\i}a Barroso, Andrzej Lenarcik and Arkadiusz P\l oski}
\date{April 2007}
\begin{document}
\maketitle
\begin{abstract}
We characterize plane curve germs nondegenerate in Kouchnirenko's
sense in terms of characteristics and intersection multiplicities
of branches.
\end{abstract}
\section{Introduction}
In this paper we consider (reduced) plane curve germs $C$,
$D$,... centered at a fixed point $O$ of a complex nonsingular
surface. Two germs $C$ and $D$ are {\it equisingular\/} if there exists a
bijection between their branches which preserves characteristic
pairs and intersection numbers. Let $(x,y)$ be a chart centered at
$O$. Then a plane curve germ has a local equation of the form $\sum
c_{\alpha,\beta}x^{\alpha}y^{\beta}=0$. Here $\sum
c_{\alpha,\beta}x^{\alpha}y^{\beta}$ is a convergent power series
without multiple factors. The {\it Newton diagram\/} $\Delta_{x,y}(C)$ is
defined to be the convex hull of the union of quadrants
$(\alpha,\beta)+(\mathbb R_+)^2$, $c_{\alpha,\beta}\neq 0$.
Recall that the {\it Newton boundary\/} $\partial\Delta_{x,y}(C)$ is the union
of the compact faces of $\Delta_{x,y}(C)$.
A germ $C$ is called {\it non-degenerate\/} with respect to
the chart $(x,y)$ if the coefficients $c_{\alpha,\beta}$ where
$(\alpha,\beta)$ runs over integral points lying on the
faces of $\Delta_{x,y}(C)$ are {\em generic} (see
Preliminaries to this Note for the precise definition). It is
well-known that the equisingularity class of a germ $C$
{\it non-degenerate\/} with respect to $(x,y)$ depends only on the Newton
polygon formed by the faces of $\Delta_{x,y}(C)$: if
$(r_1,s_1),(r_2,s_2),\ldots,(r_k,s_k)$ are subsequent vertices of
$\partial\Delta_{x,y}(C)$ then the germs $C$ and $C'$ with local
equation $x^{r_1}y^{s_1}+\cdots+ x^{r_k}y^{s_k}=0$ are equisingular.
Our aim is to give an explicit description of the non-degenerate
plane curve germs in terms of characteristic pairs and intersection
numbers of branches. In particular we show that if two germs $C$ and
$D$ are equisingular then $C$ is non-degenerate if and only if $D$
is non-degenerate. The proof of our result is based on a refined
version of Kouchnirenko's formula for the Milnor number and on the
concept of contact exponent.

\section{Preliminaries}

Let $\mathbb R_+=\{x\in\bR:\,x\geq 0\}$. For any subsets $A,B$ of the quarter
$\bR_+^2$ we consider the arithmetical sum
$A+B=\{a+b:\,a\in A\mbox{ and }b\in B\}$. If $S\subset\bN^2$ then
$\bDelta(S)$ is the convex hull of the set $S+\bR_+^2$. The subset
$\bDelta$ of $\bR_+^2$ is a {\it Newton diagram\/} if $\bDelta=\bDelta(S)$
for a set $S\subset\bN^2$ (see \cite{BK},\cite{Kouchnirenko}).
According to Teissier we put
$\{\Teis{a}{b}\}=\bDelta(S)$ if $S=\{(a,0),(0,b)\}$,
$\{\Teis{a}{\infty}\}=(a,0)+\bR_+^2$ and
$\{\Teis{\infty}{b}\}=(0,b)+\bR_+^2$ for any $a,b>0$
and call such diagrams {\it elementary Newton diagrams\/}.
The Newton diagrams form the semigroup $\cN$ with respect to the
arithmetical sum. The elementary Newton diagrams generate $\cN$.
If $\bDelta=\sum_{i=1}^r\{\Teis{a_i}{b_i}\}$ then $a_i/b_i$ are
the inclinations of edges of the diagram $\bDelta$
(by convention $\frac{a}{\infty}=0$ and
$\frac{\infty}{b}=\infty$ for $a,b>0$).
We put also $a+\infty=\infty$, $a\cdot\infty=\infty$,
$\inf\{a,\infty\}=a$ if $a>0$ and $0\cdot\infty=0$.

\medskip

\noindent {\em Minkowski's area} $[\Delta,\Delta']\in\bN\cup\{\infty\}$
of two Newton
diagrams $\Delta,\Delta'$ is uniquely determined by the following conditions

\medskip
\noindent $(m_1)$ $[\Delta_1+\Delta_2,\Delta']=[\Delta_1,\Delta']+[\Delta_2,\Delta']$,

\noindent $(m_2)$ $[\Delta,\Delta']=[\Delta',\Delta]$,

\noindent $(m_3)$ $[\{\Teis{a}{b}\},\{\Teis{a'}{b'}\}]=\hbox{\rm
inf}\{ab',a'b\}.$

\bigskip

\noindent We define the {\it Newton number\/} $\nu(\Delta)\in\bN\cup\{\infty\}$
by the following properties:

\medskip

\noindent $(\nu_1)$ $\nu(\sum_{i=1}^k \Delta_i)=\sum_{i=1}^k\nu(\Delta_i)+2\;
\sum_{1\leq i <j \leq k}[\Delta_i,\Delta_j]-k+1$,

\smallskip

\noindent $(\nu_2)$ $\nu(\{\Teis{a}{b}\})=(a-1)(b-1)$,
$\nu(\{\Teis{1}{\infty}\})=\nu(\{\Teis{\infty}{1}\})=0.$

\bigskip

\noindent
A diagram $\Delta$ is {\it convenient\/} (resp. {\it nearly convenient\/})
if $\Delta$ intersects both axes (resp. if the distances of $\Delta$
to the axes are$\mbox{}\leq 1$). Note that $\Delta$ is nearly convenient if
and only if $\nu(\Delta)\neq\infty$.
Fix a complex nonsingular surface i.e. a complex holomorphic
variety of dimension 2. In all this paper we consider {\it reduced\/}
plane curve
germs $C,D,\dots$ centered at a fixed point $O$ of this surface.
We denote by $(C,D)$ the {\it intersection multiplicity\/} of $C$ and $D$
and by $m(C)$ the {\it multiplicity\/} of $C$. We have
$(C,D)\geq m(C)m(D)$; if $(C,D)=m(C)m(D)$ then we say that $C$ and $D$
{\it intersect transversally\/}.
Let $(x,y)$ be a chart centered at $O$. Then a plane curve germ $C$ has a local
equation $f(x,y)=\sum c_{\alpha\beta}x^\alpha y^\beta\in\bC\{x,y\}$
without multiple factors. We put $\bDelta_{x,y}(C)=\bDelta(S)$ where
$S=\{(\alpha,\beta)\in\bN^2:\,c_{\alpha\beta}\neq 0\}$.
Clearly $\bDelta_{x,y}(C)$ depends on $C$ and $(x,y)$. We have two fundamental properties
of Newton diagrams:
\begin{itemize}
\item[($N_1$)] If $(C_i)$ is a finite family of plane curve germs such that
$C_i$ and $C_j$ ($i\neq j$) have no common irreducible component, then
$$\bDelta_{x,y}\left(\bigcup_iC_i\right)=\sum_i\bDelta_{x,y}(C_i)\;.$$
\item[($N_2$)] If $C$ is an irreducible germ (a branch) then
$$
\bDelta_{x,y}(C)=\left\{\Teisssr{(C,y=0)}{(C,x=0)}{10}{5}\right\}\;.
$$
\end{itemize}
\noindent For the proof we refer the reader to \cite{BK}, pp. 634--640.

\medskip

\noindent The topological boundary of $\Delta_{x,y}(C)$ is the union of two
half-lines and of a finite number of compact segments (faces). For any face
$S$ of $\Delta_{x,y}(C)$ we let
$f_S(x,y)=\sum_{(\alpha,\beta)\in S}c_{\alpha,\beta}x^{\alpha}y^{\beta}$.
Then $C$ is {\it non-degenerate\/} with respect to the chart $(x,y)$ if for all faces
$S$ of $\Delta_{x,y}(C)$ the system
$$
\frac{\partial f_S}{\partial x}(x,y)=\frac{\partial f_S}{\partial y}(x,y)=0
$$
has no solutions in $\mathbb C^*\times \mathbb C^*$.
We say that the germ $C$ is {\it non-degenerate\/} if there exists a chart $(x,y)$
such that $C$ is non-degenerate with respect to $(x,y)$.

\medskip

\noindent For any reduced plane curve germs $C$ and $D$ with
irreducible components $(C_i)$ and $(D_j)$ we put $d(C,D)=
\inf_{i,j}\{(C_i,D_j)/(m(C_i)m(D_j))\}$ and call $d(C,D)$ the {\it
order of contact\/} of germs $C$ and $D$. We have for any $C,D$ and
$E$:
\begin{itemize}
\item[($d_1$)] $d(C,D)=\infty$ if and only if $C=D$ is a branch,
\item[($d_2$)] $d(C,D)=d(D,C)$,
\item[($d_3$)] $d(C,D)\geq\inf\{d(C,E),d(E,D)\}$.
\end{itemize}
The proof of ($d_3$) is given in \cite{Chadzynski-Ploski} for the
case of irreducible $C,D,E$ which implies the general case. Condition ($d_3$)
is equivalent to the following: at least two of three numbers
$d(C,D)$, $d(C,E)$, $d(E,D)$ are equal and the third is not smaller
than the other two.
\noindent For each germ $C$ we define
$$
  d(C)=\sup\{d(C,L):\,L\mbox{ runs over all smooth branches}\}
$$
and call $d(C)$ the {\it contact exponent\/} of $C$ (see \cite{H}, Definition~1.5
where the term characteristic exponent is used). Using ($d_3$) we check that
$d(C)\leq d(C,C)$.
\begin{itemize}
\item[($d_4$)] For every finite family $(C^i)$ of plane curve germes we have
$$d(\bigcup_iC^i)=\inf\{\inf_id(C^i),\inf_{i,j}d(C^i,C^j)\}\;.$$
\end{itemize}
The proof of ($d_4$) is given in~\cite{Garcia-Lenarcik-Ploski}
(see Proposition~2.6).
We say that a smooth germ $L$ has {\it maximal contact\/} with $C$
if $d(C,L)=d(C)$. Note that $d(C)=\infty$ if and only if $C$ is a smooth
branch. If $C$ is singular then $d(C)$ is a rational number and
there exists a smooth branch $L$ which has maximal contact with $C$
(see \cite{H}, \cite{BK}).

\section{Results}
\noindent Let $C$ be a plane curve germ. A finite family of germs
$(C^{(i)})_i$ is called a {\em decomposition} of $C$ if $C=\cup_i
C^{(i)}$ and $C^{(i)},C^{(i_1)}$ $(i\neq i_1)$ have no common branch.
The following definition is basic for us.

\begin{defi}\label{Ngerm}
A plane curve $C$ is a Newton's germ (shortly N-germ) if there
exists a decomposition $(C^{(i)})_{1\leq i\leq s}$ of $C$ such that
\begin{enumerate}
\item $1\leq d(C^{(1)})<\ldots<d(C^{(s)})\leq\infty$.
\item Let $(C^{(i)}_{j})_j$ be branches of $C^{(i)}$. Then
\begin{enumerate}
\item if $d(C^{(i)})\in\bN\cup\{\infty\}$ then the branches $(C^{(i)}_{j})_j$
are smooth,
\item if $d(C^{(i)})\not\in\bN\cup\{\infty\}$ then there exists
a pair of coprime integers $(a_i,b_i)$ such that each branch
$C^{(i)}_{j}$ has exactly one characteristic pair $(a_i,b_i)$.
Moreover $d(C_j^{(i)})=d(C^{(i)})$ for all $j$.
\end{enumerate}
\item If $C^{(i)}_{l}\neq C^{(i_1)}_{k}$ then $d(C^{(i)}_{l},C^{(i_1)}_{k})=
\hbox{\rm inf}\{d(C^{(i)}),d(C^{(i_1)})\}$.
\end{enumerate}
\end{defi}

\medskip

\noindent A branch is a Newton's germ if it is smooth or has exactly
one characteristic pair. Let $C$ be a Newton's germ. The decomposition
$\{C^{(i)}\}$ satisfying (1), (2) and (3) is not unique. Take for
example a germ $C$ that has all $r>2$ branches smooth intersecting
with multiplicity $d>0$. Then for any branch $L$ of $C$ we may put
$C^{(1)}=C\setminus\{L\}$ and $C^{(2)}=\{L\}$ (or simply $C^{(1)}=C$).
If $C$ and $D$ are equisingular germs then $C$ is a $N$-germ
if and only if $D$ is a $N$-germ.

\medskip

\noindent Our main result is
\begin{teorema}\label{main-theorem}
Let $C$ be a plane curve germ. Then the following two conditions are equivalent
\begin{enumerate}
\item The germ $C$ is non-degenerate with respect to a chart $(x,y)$ such that $C$
and $\{x=0\}$ intersect transversally,
\item $C$ is a Newton's germ.
\end{enumerate}
\end{teorema}

\noindent The proof of Theorem \ref{main-theorem} we give in Section
\ref{demostracion} of this paper. Let us note here
\begin{coro} If the germ $C$ is unitangent then $C$ is non-degenerate
if and only if $C$ is a $N$-germ.
\end{coro}
Every germ $C$ has the {\em tangential decomposition\/}
$(\tilde{C}^i)_{i=1,\ldots,t}$ such that
\begin{enumerate}
\item $\tilde{C}^i$ are unitangent, that is for every two branches
$\tilde{C}^i_j$, $\tilde{C}^i_k$ of $\tilde{C}^i$ one has
$d(\tilde{C}^i_j,\tilde{C}^i_k)>1$.
\item $d(\tilde{C}^i,\tilde{C}^{i_1})=1$ for $i\neq i_1$.
\end{enumerate}

\medskip
\noindent We call $(\tilde{C}^i)_i$ tangential components of $C$.
Note that $t(C)=t$  (the number of tangential components) is an
invariant of equisingularity.

\begin{teorema}\label{thm34}
If $(\tilde{C}^i)_{i=1,\ldots,t}$ is the tangential decomposition of
the germ $C$ then the following two conditions are equivalent
\begin{enumerate}
\item  The germ $C$ is non-degenerate.
\item All tangential components $\tilde{C}^i$ of $C$ are N-germs and at
least $t(C)-2$ of them are smooth.
\end{enumerate}
\end{teorema}
Using Theorem~\ref{thm34} we get
\begin{coro}\label{cor35}
Let $C$ and $D$ be equisingular plane curve germs. Then $C$ is non-degenerate
if and only if $D$ is non-degenerate.
\end{coro}

\section{Kouchnirenko's theorem for plane curve singularities}

Let $\mu(C)$ be the {\it Milnor number\/} of a reduced germ $C$.
By definition
$\mu(C)=\mbox{dim}\,\bC\{x,y\}/(\frac{\partial f}{\partial x},%
\frac{\partial f}{\partial y})$
where $f=0$ is an equation without multiple factors of $C$.
The following properties are well-known (see for example~\cite{Ploski2}).

\begin{itemize}
\item[($\mu_1$)] $\mu(C)=0$ if and only if $C$ is a smooth branch.
\item[($\mu_2$)] If $C$ is a branch with the first characteristic
pair $(a,b)$ then $\mu(C)\geq (a-1)(b-1)$. We have
$\mu(C)=(a-1)(b-1)$ if and only if $(a,b)$ is the unique
characteristic pair of $C$.
\item[($\mu_3$)] If $(C^{(i)})_{i=1,\ldots,k}$ is a decomposition of
$C$ then $$\mu(C)=\sum_{i=1}^k \mu(C^{(i)})+2\sum_{1\leq i<j\leq
k}(C^{(i)},C^{(j)})-k+1.$$
\end{itemize}

Now we can give a refined version of Kouchnirenko's theorem
in two dimensions.
\begin{teorema}\label{Kouchnirenko}
Let $C$ be a reduced plane curve germ. Fix a
chart $(x,y)$. Then $\mu(C)\geq \nu(\Delta_{x,y}(C))$ with equality
if and only if $C$ is non-degenerate with respect to $(x,y)$.
\end{teorema}
Proof. Let $f=0$, $f\in\bC\{x,y\}$ be the local equation without
multiple factors of the germ $C$. To abbreviate the notation we put
$\mu(f)=\mu(C)$ and $\Delta(f)=\Delta_{x,y}(C)$.
If $f=x^ay^b\varepsilon(x,y)$ in $\bC\{x,y\}$ with $\varepsilon(0,0)\neq 0$
then the theorem is obvious. Then we can write
$f=x^ay^bf_1$ in $\bC\{x,y\}$ where $a,b\in\{0,1\}$ and
$f_1\in\bC\{x,y\}$ is a convenient power series. A simple
calculation based on properties ($\mu_2$), ($\mu_3$) and ($\nu_1$),
($\nu_2$) shows that $\mu(f)-\nu(\Delta(f))=\mu(f_1)-\nu(\Delta(f_1))$.
Moreover $f$ is non-degenerate if and only if if $f_1$ is non-degenerate
and the theorem reduces to the case of convenient power series
which is proved in~\cite{P} (Theorem~1.1).
\begin{nota} The implication $\mu(C)=\nu(\Delta_{x,y}(C))\Rightarrow C$
is non-degenerate is not true for the hypersurfaces with isolated singularity
{\rm(see~\cite{Kouchnirenko}, Remarque 1.21)}.
\end{nota}

\begin{coro}\label{cor43}
For any reduced germ $C$ we have $\mu(C)\geq (m(C)-1)^2$. The
equality holds if and only if $C$ is an ordinary singularity, i.e.
such that $t(C)=m(C)$.
\end{coro}

\noindent{\bf Proof.} Use Theorem  \ref{Kouchnirenko} in generic
coordinates.

\section{Proof of Theorem \ref{main-theorem}}
\label{demostracion}

\medskip
\noindent We start with the implication (1)$\Rightarrow$(2).
Let $C$ be a plane curve germ and let $(x,y)$ be a chart
such that $\{x=0\}$ and $C$ intersect transversally. The following
is well-known (\cite{Oka2}, Proposition 4.7).

\medskip
\begin{lema}
There exists a decomposition $(C^{(i)})_{i=1,\ldots,s}$ of $C$ such
that
\begin{enumerate}
\item $\Delta_{x,y}(C^{(i)})=\left\{\Teisssr{(C^{(i)},y=0)}{m(C^{(i)})}%
{12}{6}\right\}\;.$
\item Let $d_i=\frac{(C^{(i)},y=0)}{m(C^{(i)})}$. Then $1\leq
d_1<\cdots<d_s\leq\infty$ and $d_s=\infty$ if and only if
$C^{(s)}=\{y=0\}$.
\item Let $n_i=m(C^{(i)})$ and $m_i=n_id_i=(C^{(i)},y=0)$. Suppose that $C$
is non-degenerate with respect to the chart $(x,y)$. Then $C^{(i)}$
has $r_i=\hbox{\rm g.c.d.}(n_i,m_i)$ branches $C^{(i)}_j:
y^{n_i/r_i}-a_{ij}x^{m_i/r_i}+\cdots=0$ {\rm(}$j=1,\ldots,r_i$ and
$a_{ij}\neq a_{ij'}$, if $j\neq j'${\rm)}.
\end{enumerate}
\end{lema}

\medskip
\noindent Using the above lemma we prove that any germ $C$
non-degenerate with respect to $(x,y)$ is a $N$-germ.
From ($d_4$) we get $d(C^{(i)})=d_i$.
Clearly all branches $C^{(i)}_j$ have exactly one
characteristic pair $(\frac{n_i}{r_i},\frac{m_i}{r_i})$ or are
smooth. A simple calculation shows that
$$
d(C^{(i)}_j,C^{(i_1)}_{j_1})=\frac{(C^{(i)}_j,C^{(i_1)}_{j_1})}
{m(C^{(i)}_j)m(C^{(i_1)}_{j_1})}=\hbox{\rm inf}\{d_i,d_{i_1}\}\;.
$$
To prove the implication (2)$\Rightarrow$(1) we need some auxilary Lemmas.
\begin{lema}\label{lm52}
Let $C$ be a plane curve germ which all branches $C_i$ $(i=1,\dots,s)$
are smooth. Then there exists a smooth germ $L$ such that
$(C_i,L)=d(C)$ for $i=1,\dots,s$.
\end{lema}
Proof. If $d(C)=\infty$ then $C$ is smooth and we take $L=C$.
If $d(C)=1$ then we take a smooth germ $L$ such that $C$ and $L$ are
transversal. Let $k=d(C)$ and suppose that $1<k<\infty$.
By formula ($d_4$) we get $\inf\{(C_i,C_j):\,i,j=1,\dots,s\}=k$.
We may assume that $(C_1,C_2)=\dots=(C_1,C_r)=k$ and $(C_1,C_j)>k$
for $j>r$ for an index $r$, $1\leq r\leq s$.
There is a system of coordinates $(x,y)$ such that $C_j$ ($j=1,\dots,r$)
have equations $y=c_jx^k+\dots$. It suffices to take $L:\,y-cx^k=0$
where $c\neq c_j$ for $j=1,\dots,r$.
\begin{lema}\label{lm53}
Suppose that $C$ is a $N$-germ and let $(C^{(i)})_{1\leq i\leq s}$ be
a decomposition of $C$ such as in Definition~\ref{Ngerm}.
Then there is a smooth germ $L$ such that $d(C_j^{(i)},L)=d(C^{(i)})$
for all $j$.
\end{lema}
Proof. Step 1. There is a smooth germ $L$ such that
$d(C_j^{(s)},L)=d(C^{(s)})$ for all $j$. If $d(C^{(s)})\in\bN\cup\{\infty\}$
then the existence of $L$ follows from Lemma~\ref{lm52}.
If $d(C^{(s)})\notin\bN\cup\{\infty\}$ then all components $C^{(s)}_j$
have the same characteristic pair $(a_s,b_s)$. Fix a component $C^{(s)}_{j_0}$
and let $L$ be a smooth germ such that
$d(C^{(s)}_{j_0},L)=d(C^{(s)}_{j_0})=d(C^{(s)})$.
Let $j_1\neq j_0$. Then
$d(C^{(s)}_{j_1},L)\geq\inf\{d(C^{(s)}_{j_1},C^{(s)}_{j_0}),%
d(C^{(s)}_{j_0},L)\}=d(C^{(s)})$. On the other hand
$d(C^{(s)}_{j_1},L)\leq d(C^{(s)}_{j_1})=d(C^{(s)})$ and we get
$d(C^{(s)}_{j_1},L)=d(C^{(s)})$.

\noindent Step 2. Let $L$ be a smooth germ such that
$d(C^{(s)}_j,L)=d(C^{(s)})$ for all $j$. We will check that
$d(C_j^{(i)},L)=d(C^{(i)})$ for all $i$ and $j$. To this purpose
fix $i<s$. Let $C^{(s)}_{j_0}$ be a component of $C^{(s)}$. Then
$d(C^{(i)}_j,C^{(s)}_{j_0})=\inf\{d(C^{(i)}),d(C^{(s)})\}=d(C^{(i)})$.
By ($d_3$) we get $d(C^{(i)}_j,L)\geq\inf\{d(C^{(i)}_j,C^{(s)}_{j_0}),%
d(C^{(s)}_{j_0},L)\}=\inf\{d(C^{(i)}),d(C^{(s)})\}=d(C^{(i)})$.
On the other hand $d(C^{(i)}_j,L)\leq d(C^{(i)}_j)=d(C^{(i)})$
and we are done.
\begin{nota}\label{rem54}
In notation of the above lemma we have
$(C^{(i)},L)=m(C^{(i)})d(C^{(i)})$ for $i=1,\dots,s$.
\end{nota}
Indeed, if $C^{(i)}_j$ are branches of $C^{(i)}$ then
{\small
$$
(C^{(i)},L)=\sum_j(C^{(i)}_j,L)=\sum_j m(C^{(i)}_j)d(C^{(i)}_j,L)%
=\sum_j m(C^{(i)}_j)d(C^{(i)})=m(C^{(i)})d(C^{(i)})\;.
$$
} 
\begin{lema}\label{lm55}
Let $C$ be a $N$-germ and let $(C^{(i)})_{1\leq i\leq s}$ be a decomposition
of $C$ such as in Definition~\ref{Ngerm}. Then
{\small $$
\mu(C)=\sum_i(m(C^{(i)})-1)(m(C^{(i)})d(C^{(i)})-1)+2\sum_{i<j}%
m(C^{(i)})m(C^{(j)})\inf\{d(C^{(i)}),d(C^{(j)}\}-s+1\;.
$$ } 
\end{lema}
Proof. Use properties ($\mu_1$),($\mu_2$) and ($\mu_3$) of the Milnor number.

\medskip
\noindent To prove implication (2)$\Rightarrow$(1) of Theorem~\ref{main-theorem}
suppose that $C$ is
a $N$-germ and let $(C^{(i)})_{i=1,\dots,s}$ be a decomposition of $C$
such as in Definition~\ref{Ngerm}.
Let $L$ be a smooth branch such that $(C^{(i)},L)=m(C^{(i)})d(C^{(i)})$
for $i=1,\dots,s$ (such a branch exists by Lemma~\ref{lm53} and Remark~\ref{rem54}).
Take a system of coordinates such that $\{x=0\}$ and $C$ are transversal
and $L=\{y=0\}$. Then we get
$$
\Delta_{x,y}(C)=\sum_{i=1}^s\Delta_{x,y}(C^{(i)})=
\sum_{i=1}^s\left\{\Teisssr{(C^{(i)},\{y=0\})}{m(C^{(i)})}{14}{7}\right\}=
\sum_{i=1}^s\left\{\Teisssr{m(C^{(i)})d(C^{(i)})}{m(C^{(i)})}{14}{7}\right\}
$$
and consequently
\begin{eqnarray*}
\nu(\Delta_{x,y}(C)) & = & \sum_{i=1}^s(m(C^{(i)})-1)(m(C^{(i)})d(C^{(i)})-1)\\
& &\quad\mbox{}+2\sum_{1\leq i<j\leq s}m(C^{(i)})m(C^{(j)})\inf\{d(C^{(i)}),d(C^{(j)})\}-s+1\\
& = & \mu(C)
\end{eqnarray*}
by Lemma~\ref{lm55}. Therefore $\mu(C)=\nu(\Delta_{x,y}(C))$ and $C$ is
non-degenerate with respect to $(x,y)$ by Theorem~\ref{Kouchnirenko}.

\section{Proof of Theorem~\ref{thm34}}
The Newton number $\nu(C)$ of the plane curve germ $C$ is defined to
be $\nu(C)=\hbox{\rm sup}\{\nu(\Delta_{x,y}(C))\;:\;(x,y)
\;\;\hbox{\rm runs over all charts centered at $O$}\}$.

\medskip
\noindent Using Theorem \ref{Kouchnirenko} we get
\begin{lema}\label{lm61}
A plane curve germ $C$ is non-degenerate if and only if
$\nu(C)=\mu(C)$.
\end{lema}
The proposition below shows that we can reduce the computation
of the Newton number to the case of unitangent germs.
\begin{prop}\label{prop62}
\label{general} If $C=\bigcup_{k=1}^t \tilde{C}^k$ $(t>1)$ where
$\{\tilde{C}^k\}_k$ are unitangent germs such that
$(\tilde{C}^k,\tilde{C}^l)=m(\tilde{C}^k)m(\tilde{C}^l)$
for $k\neq l$ then
$$\nu(C)-(m(C)-1)^2=\hbox{\rm max}_{1\leq k<l\leq t}
\{(\nu(\tilde{C}^k)-(m(\tilde{C}^k)-1)^2)+(\nu(\tilde{C}^l)-(m(\tilde{C}^l)-1)^2)\}.$$
\end{prop}
Proof. Let $\tilde{n}_k=m(\tilde{C}^k)$. Suppose
that $\{x=0\}$ and $\{y=0\}$ are tangent to $C$. Then there are two
tangencial components $\tilde{C}^{k_1}$ and $\tilde{C}^{k_2}$ such
that $\{x=0\}$ is tangent to $\tilde{C}^{k_1}$ and $\{y=0\}$ is
tangent to $\tilde{C}^{k_2}$. Now, we have

\begin{eqnarray*}
\nu(\Delta_{x,y}(C))&=&\nu(\sum_{k=1}^t \Delta_{x,y}(\tilde{C}^k))
=\nu(\Delta_{x,y}(\tilde{C}^{k_1}))+\nu(\Delta_{x,y}(\tilde{C}^{k_2}))\\
&+&\sum_{k\neq k_1,k_2}\nu(\Delta_{x,y}(\tilde{C}^{k}))
+2\sum_{1\leq k<l\leq t}
\left[\Delta_{x,y}(\tilde{C}^{k}),\Delta_{x,y}(\tilde{C}^{l})\right]-t+1\\
 &=&\nu(\Delta_{x,y}(\tilde{C}^{k_1}))+\nu(\Delta_{x,y}(\tilde{C}^{k_2}))+
\sum_{k\neq k_1,k_2}(\tilde{n}_k-1)^2+2\sum_{1\leq k< l \leq t}\tilde{n}_k \tilde{n}_l-t+1\\
&=&\nu(\Delta_{x,y}(\tilde{C}^{k_1}))-(\tilde{n}_{k_1}-1)^2+\nu(\Delta_{x,y}(\tilde{C}^{k_2}))-(\tilde{n}_{k_2}-1)^2+(m(C)-1))^2.
\end{eqnarray*}
The germs $\tilde{C}^{k_1}$ and $\tilde{C}^{k_2}$ are
unitangent and transversal. Thus it is easy to see that there exists
a chart $(x_1,y_1)$ such that
$\nu(\Delta_{x_1,y_1}(\tilde{C}^k))=\nu(\tilde{C}^k)$ for
$k=k_1,k_2.$

\noindent If $\{x=0\}$ (or $\{y=0\}$) and $C$ are transversal then
there exists a $k\in \{1,\ldots,t\}$ such that
$\nu(\Delta_{x,y}(C))=\nu(\Delta_{x,y}(\tilde{C}^{k}))-(\tilde{n}_{k}-1)^2+(m(C)-1))^2$
and the proposition follows from the previous considerations.

Now we can pass to the proof of Theorem~\ref{thm34}. If $t(C)=1$ then $C$ is
non-degenerate with respect to a chart $(x,y)$ such that $C$ and
$\{x=0\}$ intersect transversally and Theorem~\ref{thm34} follows from
Theorem~\ref{main-theorem}. If $t(C)>1$ then
by Proposition~\ref{prop62} there are indices
$k_1<k_2$ such that
\begin{itemize}
\item[($\alpha$)]
$\nu(C)-(m(C)-1)^2=\nu(\tilde{C}^{k_1})-(m(\tilde{C}^{k_1})-1)^2+
                  \nu(\tilde{C}^{k_2})-(m(\tilde{C}^{k_2})-1)^2\;.$
\end{itemize}
On the other hand from basic properties of the Milnor number we get
\begin{itemize}
\item[($\beta$)]
$\mu(C)-(m(C)-1)^2=\sum_k(\mu(\tilde{C}^k)-(m(\tilde{C}^k)-1)^2)\;.$
\end{itemize}
Using ($\alpha$), ($\beta$) and Lemma~\ref{lm61} we check that $C$ is non-degenerate
if and only if $\mu(\tilde{C}^{k_1})=\nu(\tilde{C}^{k_1})$,
$\mu(\tilde{C}^{k_2})=\nu(\tilde{C}^{k_2})$
and $\mu(\tilde{C}^k)=(m(\tilde{C}^k)-1)^2$ for $k\neq k_1,k_2$. Now
Theorem~\ref{thm34} follows from Lemma~\ref{lm61} and Corollary~\ref{cor43}.

\section{Concluding remark}
M.~Oka proved in~\cite{Oka1} that the Newton number like the Milnor
number is an invariant of equisingularity. Therefore the invariance
of non-degeneracy (Corollary~\ref{cor35}) follows from the equality
$\nu(C)=\mu(C)$ characterizing non-degenerate germs (Lemma~\ref{lm61}).

\section*{Acknowledgements}
The third author (A.P.) is grateful to La Laguna University
where a part of this work was prepared.

\medskip
\noindent
{\small Evelia Rosa Garc\'{\i}a Barroso\\
Departamento de Matem\'atica Fundamental\\
Facultad de Matem\'aticas, Universidad de La Laguna\\
38271 La Laguna, Tenerife, Espa\~na\\
e-mail: ergarcia@ull.es}

\medskip
\noindent {\small Andrzej Lenarcik\\
Department of Mathematics\\
Technical University \\
Al. 1000 L PP7\\
25-314 Kielce, Poland\\
e-mail: ztpal@tu.kielce.pl}

\medskip
\noindent {\small Arkadiusz P\l oski\\
Department of Mathematics\\
Technical University \\
Al. 1000 L PP7\\
25-314 Kielce, Poland\\
e-mail: matap@tu.kielce.pl}

\end{document}